\documentclass[reqno,A4,11pt,newtimes]{amsart}
\usepackage{mathrsfs}
\usepackage{amsfonts}
\usepackage{color}
\usepackage{enumerate}
\usepackage{hyperref}
\usepackage{indentfirst}
\usepackage{amsmath}
\usepackage{amssymb}
\usepackage{enumerate}

\hypersetup{hidelinks}
\usepackage{latexsym,eucal,amsfonts}
\usepackage{color,graphics,graphicx}
\def\qed{\hfill \rule{4pt}{7pt}}

\def\pf{\noindent {\it Proof.} }
\def\hrulefill{\leavevmode\leaders\hrule height 4pt\hfill\kern0pt}
\makeatletter
\def\rulefill{\leavevmode\leaders\hrule depth -3pt height 6pt\hfill\kern0pt}
\newcommand{\poq}[2]{(#1;q)_{#2}}
\newcommand{\poqq}[2]{(#1;q^2)_{#2}}

\newcommand{\be}{\begin{equation}}
\newcommand{\ee}{\end{equation}}
\newcommand{\ba}{\begin{array}}
\newcommand{\ea}{\end{array}}
\newcommand{\bmn}{\begin{eqnarray}}
\newcommand{\emn}{\end{eqnarray}}
\newcommand{\bnm}{\begin{eqnarray*}}
\newcommand{\enm}{\end{eqnarray*}}
\newcommand{\bln}{\begin{subequations}}
\newcommand{\eln}{\end{subequations}}

\newtheorem{thm}{Theorem}[section]
\newtheorem{lemma}[thm]{Lemma}
\newtheorem{cor1}[thm]{Corollary}
\newtheorem{cor}[thm]{Theorem}
\newtheorem{prop}[thm]{Proposition}
\newtheorem{exam}[thm]{Example}
\newtheorem{remark}[thm]{Remark}

\numberwithin{equation}{section}

\usepackage{amstext}
\begin{document} 

{
\title{Some $q$-Identities derived by the ordinary derivative operator}
\author{ J. Wang$^{a,\ast}$}
\author{R. Q. Ruan$^{b}$}
\author{X. R. Ma$^{c}$}
\dedicatory{$^{a,b}$ School of Mathematical Sciences, Zhejiang Normal  University, Jinhua 321004, P.R.China\\ $^{c}$Department of Mathematics, Soochow University, SuZhou 215006, P.R.China}
\thanks{*Corresponding author. E-mail address: jinwang@zjnu.edu.cn}
\thanks{}
\begin{abstract}
In this paper, we investigate applications of the ordinary derivative operator, instead of the $q$-derivative operator,  to the theory of $q$-series. As main results, many new  summation and transformation formulas   are established  which  are closely related to some well-known formulas such as the $q$-binomial theorem, Ramanujan's ${}_1\psi_1$ formula, the quintuple product identity, Gasper's $q$-Clausen product formula, and Rogers' ${}_6\phi_5$  formula, etc. Among these results is a finite form of the Rogers-Ramanujan identity and a short way to Eisenstein’s theorem on Lambert series.
\end{abstract}
\maketitle\thispagestyle{empty}
\markboth{J. Wang, R. Q. Ruan, X. R. Ma}
{Some $q$-identities derived by the ordinary derivative operator }

\section{Introduction}
\subsection{Two derivative operators}
It is well known that in classical calculus, the derivative of a function $f(x)$ is defined to be
$$
\mathbf{D}_x(f(x)):=\frac{d f}{d x}=\lim _{h \rightarrow 0} \frac{f(x+h)-f(x)}{h}.
$$
Naturally,  it can be extended to   the $q$-derivative operator $\mathbf{D}_{q,x}$ \cite[Ex. 1.12]{Gasper}
$$
\mathbf{D}_{q,x}(f(x))=\frac{f(q x)-f(x)}{q x-x}.
$$
As of today,  the $\mathbf{D}_{q,x}(f(x))$ introduced first by  H. Jackson \cite{jackson} plays a very important role in  the theory of $q$-Analysis or $q$-Calculus, as well as the Quantum Calculus. We refer the reader to T. Ernst's  paper \cite{ernst}  for a history of $q$-Calculus and to V. Kac and P. Cheung's book \cite{kac} for deep applications in Physics. In what follows,  we call $\mathbf{D}_{x}$  the ordinary differential (in short, OD) operator  and  write
\begin{align*} \mathbf{D}_{x=a}(f(x)):=\mathbf{D}_{x}(f(x))|_{x=a}.
\end{align*}
It is obvious that
\begin{align*}
\lim_{q\to 1}\mathbf{D}_{q,x}(f(x))=\mathbf{D}_x(f(x)).
\end{align*}
 Roughly speaking, in  the theory of $q$-Analysis,  $\mathbf{D}_{q,x}$ has attracted much attention and a lot of results have been found \cite{Gasper}. By contrast, the OD-operator $\mathbf{D}_x$ seems to have been  ignored and less  discussed, let along the related results. Although  $\mathbf{D}_{x}$  is the limiting case  of $\mathbf{D}_{q,x}$   as $q$ tends to 1, when applied to the same $q$-identity,  it often leads us to different  results. Regarding this, we had better  give a concrete example.
\begin{exam}Consider the basic $q$-identity
\begin{align}
\sum_{n=0}^{\infty}\frac{x^n}{(q;q)_n}=\frac{1}{\poq{x}{\infty}}\quad(|x|<1).\label{eqx}
\end{align}
By applying $\mathbf{D}_x$ and $\mathbf{D}_{q,x}$ to both sides of \eqref{eqx}, respectively, we have the following results:  \begin{description}
\item[(1) Action of $\mathbf{D}_x$]
by applying the OD-operator $\mathbf{D}_x$ to  \eqref{eqx}, we obtain
\begin{align}
\sum_{n=1}^{\infty}\frac{nx^{n-1}}{(q;q)_n}=\frac{1}{\poq{x}{\infty}}\sum_{n=0}^{\infty}\frac{q^n}{1-xq^n}.\label{cankaoct}
\end{align}
\item[(2) Action of $\mathbf{D}_{q,x}$] alternatingly, if we apply the $q$-derivative operator $\mathbf{D}_{q,x}$ to \eqref{eqx}, because
\begin{align*}
\mathbf{D}_{q,x}(x^n)=\frac{1-q^n}{1-q}x^{n-1},\end{align*} then it follows
\begin{align}
\sum_{n=1}^{\infty}\frac{1}{(q;q)_n}\frac{1-q^n}{1-q}x^{n-1}=\frac{1}{1-q}\frac{1}{\poq{x}{\infty}}.\label{cankao}
\end{align}
\end{description}
It can be checked with easy
$$\lim_{q\to 1^{-}}
\mathbf{D}_{q,x}(x^n)=x^{n-1}\lim_{q\to 1^{-}}\frac{1-q^n}{1-q}=nx^{n-1},$$
but  \eqref{cankao} is not the limitation of  \eqref{cankaoct} when $q\to 1^{-}$.
\end{exam}
It is also worth remarking that if we make some obvious simplifications
in \eqref{cankao}, it becomes \eqref{eqx} again, so that $1/(x; q)_{\infty}$ is essentially an exponential function for the $q$-derivative.

From this simple example, we feel that the OD-operator $\mathbf{D}_x$, just as the $q$-derivative operator $\mathbf{D}_{q,x}$ does,  has its own value in the study of $q$-series. In practice,  as tools for $q$-series, $\mathbf{D}_x$ is not so popular as $\mathbf{D}_{q,x}$ does. Hence, we think it is worthwhile to investigate  systematically applications of $\mathbf{D}_x$ to this field. It is the main theme of the present paper.
\subsection{Notation}
Some remarks  on the notation are necessary. Throughout this paper we will use the standard notation and terminology
for $q$-series as in G. Gasper and M. Rahman
\cite{Gasper}. As customary, the $q$-shifted factorials are
\begin{align*}
(a ; q)_\infty &:=\prod_{k=0}^\infty\left(1-a q^{k}\right),\\
(a ; q)_n &:=(a ; q)_{\infty} /\left(a q^n ; q\right)_{\infty} \quad (n\in \mathbb{Z}),
\end{align*}
where $\mathbb{Z}$ denotes the set of all integers. In general, the multi-factorial is given by
\begin{align*}
\left(a_1, a_2, \ldots, a_m ; q\right)_n:=\prod_{k=1}^m\left(a_k; q\right)_n.
\end{align*}
As usual,  the binomial coefficient $\binom{n}{k}=n!/(k!(n-k)!)$ and its $q$-analogue  is
$$\left[n \atop k\right]_q :=\frac{(q ; q)_n}{(q ; q)_k(q ; q)_{n-k}} .
$$
Usually it is referred to as the $q$-binomial coefficient or the Gaussian binomial coefficient. Following \cite{Gasper}, we introduce a basic hypergeometric series with the base $q: |q|<1$ and the argument $x$  as
\begin{align}
{ }_{r} \phi_s\left[\begin{array}{c}
a_1, a_2,\ldots, a_{r} \\
b_1,b_2, \ldots, b_s
\end{array}; q, x\right]:=\sum_{n=0}^{\infty} \frac{\left(a_1,a_2,\ldots, a_{r} ; q\right)_n}{\left(b_1, b_2,\ldots, b_s ; q\right)_n}\big((-1)^nq^{\binom{n}{2}}\big)^{s+1-r} \frac{x^n}{(q; q)_n}.\label{notation-two}
\end{align}
In the meantime, the bilateral hypergeometric series for $ x \neq 0$ is defined to be
\begin{align}
{ }_r \psi_r\left[\begin{array}{l}
a_1,a_2, \ldots, a_r \\
b_1,b_2, \ldots, b_r
\end{array}; q, x\right]:=\sum_{n=-\infty}^{\infty} \frac{\left(a_1, a_2,\ldots, a_r ; q\right)_n}{\left(b_1, b_2,\ldots, b_r ; q\right)_n} x^n.\label{notation-three}
\end{align}
\subsection{The method how to use the OD-operator $\mathbf{D}_{x}$}
Our method consists of the following facts.
\begin{prop}[Leibniz's rule]\label{generalpri-1}
For any finite integer $n\geq 0$ and two functions $F(x)$ and $G(x)$ are  differentiable $n$ times at $x$, it holds
\begin{align}\mathbf{D}_{x}^n(F(x)G(x))=\sum_{k=0}^n\binom{n}{k}\mathbf{D}_{x}^k(F(x))\mathbf{D}_{x}^{n-k}(G(x)).
\end{align}
\end{prop}
As far as $q$-series is concerned, we often deal with a given $q$-identity via  Leibiniz's rule in such way.
\begin{prop}\label{generalpri}  Suppose that there holds the following $q$-identity
\begin{align}\sum_{n\in \mathbb{Z}}T_n(x)=\poq{ax}{\infty}S(x),\label{vvvvvv}
\end{align}
where $S(q^{-m}/a)\neq 0$  for any finite integer $m\geq 0$. Then we have
\begin{align}
\sum_{n\in \mathbb{Z}}\mathbf{D}_{x=q^{-m}/a}(T_n(x))=a(-1)^{m-1}q^{-m(m-1)/2}\poq{q}{m}\poq{q}{\infty}S(q^{-m}/a).\label{vvvvvv-dual}
\end{align}
\end{prop}
\pf~It suffices to differentiate \eqref{vvvvvv} with respect to  $x$. By Leibniz's rule, it is easy to obtain
\begin{align*}
\sum_{n\in \mathbb{Z}} \mathbf{D}_x(T_n(x))=(ax;q)_\infty~ \mathbf{D}_x(S(x))+ \mathbf{D}_x((ax;q)_\infty)~S(x).
\end{align*}
Next, set  $x=q^{-m}/a$. Then
\begin{align*}
\sum_{n\in \mathbb{Z}} \mathbf{D}_{x=q^{-m}/a}(T_n(x))=a(-1)^{m-1}q^{-m(m-1)/2}\poq{q}{m}\poq{q}{\infty}S(q^{-m}/a).
\end{align*}
Note that in this procession, we have used
\begin{align*} \mathbf{D}_x((ax;q)_\infty)&= \mathbf{D}_x((1-axq^m)L(x))\\
&=(-aq^m)L(x)+(1-axq^m) \mathbf{D}_x(L(x)),\end{align*}
where $L(x)=(ax;q)_m(axq^{m+1};q)_\infty.$
\qed

As is expected,   once   $S(x)$ and $T(x)$ in \eqref{vvvvvv} are appropriately chosen,  we will be  led  to a possibly new $q$-identity \eqref{vvvvvv-dual}. Besides  Propositions \ref{generalpri-1} and \ref{generalpri}, a very crucial fact for our method is that once the OD-operator $\mathbf{D}_x$ is applied to the $q$-factorial, we have the following  computational  results. They will be frequently used in our discussions.
\begin{prop}\label{generalpri-3}\label{qiudao}
 Define $F(x):=(x ; q)_{n}$ and
 \begin{align}
F_k(x):=\prod_{i=0,i\neq k}^{n-1}(1-xq^i).\label{qiudaoeq4}
\end{align}
Then for any integer $n\geq 1$, it holds
\begin{align}
& \mathbf{D}_x(F(x))=-\sum_{k=0}^{n-1}q^k  F_k(x) ,\\
& \mathbf{D}_{x=1}(F(x))=-(q, q)_{n-1}.
\end{align}
\end{prop}
\pf All these can be verified by the usual computation, which are omitted here.
\qed

The rest of our paper is mainly concerned with applications of Propositions \ref{generalpri-1}-\ref{generalpri-3} to the theory of $q$-series.  As we will see later,  many interesting $q$-identities (stated  by theorems) are  deduced from some known formulas (quoted as lemmas)  which are listed in {\sl\bf Digital Library of Mathematical Functions}, i.e., \url{http://dlmf.nist.gov/17} or \cite[Chap.17]{dlmf}.
To the best of our knowledge,  most of these theorems are new. Even if  not new, the corresponding  proofs are novel.
\section{Main results-I: $q$-Identities from some classical formulas}
In this section, we will focus on some classical but seemingly succinct formulas  from \url{http://dlmf.nist.gov/17}, in order to pursue any new result.
\subsection{$q$-Identities based on the $q$-binomial theorem}
Let us begin with
\begin{lemma}[The $q$-binomial theorem: \url{http://dlmf.nist.gov/17.2.E35}] For any integer $n\geq 0$, it holds
\begin{align}
(x; q)_n=\sum_{k=0}^n\left[n \atop k\right]_q (-1)^kq^{\binom{k}{2}}x^k.\label{xanxan}
\end{align}
\end{lemma}
As we will see,  the action of the OD-operator $\mathbf{D}_x$ to \eqref{xanxan} gives rise to
\begin{cor}\label{fffggg}For any integers $n-1\geq m\geq  0$, it holds
\begin{align}
\sum_{k=1}^nk (-1)^kq^{\binom{k}{2}} \left[n \atop k\right]_q  q^{-mk}&=(-1)^{m+1}q^{-\binom{m+1}{2}}\poq{q}{m} \poq{q}{n-m-1}.
\end{align}
\end{cor}
\pf~ Following the aforementioned method, we differentiate both sides of \eqref{xanxan} with respect to  $x$. The result is
\begin{align}-\sum_{k=0}^{n-1} F_k(x)q^k=\sum_{k=1}^nk (-1)^kq^{\binom{k}{2}}\left[n \atop k\right]_q  x^{k-1},\label{idd13}
\end{align}
where, for $0\leq k\leq n-1$,
\begin{align*} F_k(x):=\frac{(x ; q)_n}{1-x q^k}.
\end{align*}
Then, by Proposition \ref{qiudao}, it is clear that
\begin{align*}
F_k(q^{-m})=\left\{\begin{array}{cc}
0, &0\leq k\neq m\leq n-1;\\
(-1)^m q^{- {m\choose 2}-m}\poq{q}{m} \poq{q}{n-1-m},&k=m.
\end{array}
\right.
\end{align*}
A direct substitution of these results in \eqref{idd13} with $x=q^{-m}$ yields \begin{align*}
\sum_{k=1}^nk(-1)^kq^{\binom{k}{2}} \left[n \atop k\right]_q q^{-m(k-1)}= (-1)^{m+1} q^{-{m\choose 2}}\poq{q}{m} \poq{q}{n-1-m} .
\end{align*}
Identity \eqref{fffggg} is thereby proved.
\qed

Using the same argument as above, we can show
\begin{cor}For any integer $n\geq 1$,  it holds
\begin{align}
\sum_{k=1}^nk\left[
n \atop
k
\right]_qq^{\binom{k}{2}} (-ax)^{k-1}
 =(a x; q)_n \sum_{k=0}^{n-1} \frac{q^k }{1-a x q^k}.\label{eq244}
\end{align}
\end{cor}
\pf Performing as in Theorem \ref{fffggg}, we still apply $\mathbf{D}_{x}$ to
another form of  the $q$-binomial theorem
\begin{align*}
 (a x ; q)_ n=\sum_{k=0}^n\left[
n \atop
k
\right]_q (-1)^kq^{\binom{k}{2}}(a x)^k.
\end{align*}
It is easy to compute
\begin{align*} \mathbf{D}_{x}\bigg(\sum_{k=0}^n\left[
n \atop
k
\right]_q (-1)^kq^{\binom{k}{2}}(a x)^k\bigg)=\sum_{k=1}^n\left[
n \atop
k
\right]_q (-1)^kq^{\binom{k}{2}} k a^k x^{k-1}
\end{align*}
while
\begin{align*}
\mathbf{D}_{x}\big((a x, q)_n\big)=\sum_{k=0}^{n-1}\frac{(a x; q)_n}{1-a x q^k}\mathbf{D}_{x}\left(1-a x q^k\right)
=(a x ; q)_n \sum_{k=0}^{n-1} \frac{-a q^k}{1-a x q^k}.
\end{align*}
Thus \eqref{eq244} follows.
\qed

It is of interest to see that  by taking $(x, a)=(1,q)$ and then letting $n\to\infty$, we recover Eisenstein's theorem for Lambert series.
\begin{cor1}[Eisenstein's theorem: \mbox{\rm \cite[Thm.39]{johnson}}] Let $d(n)$ be the number of divisors of $n$. Then
\begin{align}
\sum_{k=0}^\infty k\frac{(-1)^{k}q^{k(k+1)/2}}{\poq{q}{k}}
 =-(q; q)_\infty\sum_{k=1}^{\infty} d(k) q^k.
\end{align}
\end{cor1}
\pf It  is the direct consequence of \eqref{eq244} via the use of  Lambert's theorem (\cite[Thm.36]{johnson}):
\[\sum_{k=1}^{\infty} \frac{q^{k} }{1-q^k}=\sum_{k=1}^{\infty} d(k) q^k .\]
\qed

Continuing in this way, we now consider the  second-order (even more) derivatives   in $x$ at $x=1$. Hence, we achieve  the following $q$-identity related to the Lambert series.
\begin{cor}For any integer $n\geq 1$, it holds
\begin{align}\sum_{k=0}^n\binom{k}{2}\left[n \atop k\right]_q (-1)^kq^{\binom{k}{2}}=
 (q ; q)_{n-1} \sum_{k=1}^{n-1} \frac{q^k}{1-q^k}.\label{idd14}
\end{align}
In particular,
\begin{align}\sum_{k=0}^{\infty}\binom{k}{2}\frac{(-1)^kq^{\binom{k}{2}}}{(q;q)_k}=
 (q ; q)_{\infty} \sum_{k=1}^{\infty} \frac{q^k}{1-q^k}.\label{idd14-14}
\end{align}
\end{cor}
\pf~To show \eqref{idd14}, we only need to compute the second-order derivatives of both sides of the $q$-binomial theorem \eqref{xanxan}
 and then set $x=1$ in the resulting identity. To do that, by Leibniz's rule, it is easy to find
\begin{align*}
\mathbf{D}_x^2((x ; q)_n)&=\mathbf{D}_x\{-(xq;q)_{n-1}+(1-x)~\mathbf{D}_x((xq;q)_{n-1})\}\\
&=-2\mathbf{D}_x\left((xq;q)_{n-1}\right)+(1-x)\mathbf{D}_x^2\left((xq;q)_{n-1}\right).
\end{align*}
Consequently, by virtue of Proposition \ref{qiudao}, we have
\begin{align*}
 \mathbf{D}_{x=1}^2\left((x ; q)_n\right)&=-2\mathbf{D}_{x=1}\left((xq;q)_{n-1}\right)
\\
&=2(q ; q)_{n-1} \sum_{k=1}^{n-1} \frac{q^k}{1-q^k}.
\end{align*}
It completes the proof of \eqref{idd14}. Identity \eqref{idd14-14} is the limiting case of \eqref{idd14} for $n\to\infty$.
\qed

\subsection{$q$-Identities  based on the $_1\phi_0$ formula}
 It is known that \begin{align}
{ }_1 \phi_0\left[\begin{array}{l}
a \\
-
\end{array} ;q, x\right]=\frac{(ax; q)_{\infty}}{(x; q)_{\infty}}.\label{eq288}
\end{align}
 See \url{https://dlmf.nist.gov/17.5.E2} for detail.  Treated as  an analytic function in  variable  $x$, it is easy to find  its  first-order derivative, which leading us to
\begin{cor} For any finite integer $m\geq 0$ and complex number $
a$ with $|1/aq^m|<1$, we have
\begin{align}
 \sum_{n=0}^{\infty}n \frac{(a ; q)_n}{(q ; q)_n}\left(\frac{1}{a q^m}\right)^n&=(-1)^{m+1}q^{-m(m+1)/2}(q ; q)_m \frac{(q; q)_{\infty}}{\left(1/aq^m ; q\right)_\infty}.\label{eq2.8}
\end{align}
\end{cor}
\pf It follows from  applying $\mathbf{D}_{x=q^{-m}/a}$ to both sides of  \eqref{eq288} directly.
\qed
\subsection{$q$-Identities based on the ${}_1\phi_1$ formula}
Among the various formulas,  Cauchy's formula is more basic.
\begin{lemma}[cf.  \url{http://dlmf.nist.gov/17.5.E5}]
\begin{align}
{ }_1 \phi_1\left[\begin{array}{l}
a \\
c
\end{array} ;q,\frac{c}{a}\right]=\frac{(c / a ; q)_{\infty}}{(c ; q)_{\infty}}.\label{CCCC}
\end{align}
\end{lemma}
Upon applying the OD-operator  to \eqref{CCCC}, we at once obtain
\begin{cor}
For complex number  $|c|>1$, it holds
\begin{align}
\sum_{n=1}^\infty\frac{(-1)^{n-1}q^{\binom{n}{2}}}{\poq{q}{n}}\bigg\{\sum_{k=0}^{n-1}\frac{1}{1-cq^k}\bigg\}
=\frac{(q; q)_{\infty}}{(c ; q)_{\infty}}.\label{eq8}
\end{align}
\end{cor}
\pf~It suffices to take the derivative on both sides of \eqref{CCCC} with respect to $a$ at $a=c$. Then we have
 \begin{align*}
&-\sum_{n=1}^\infty(-1)^nq^{\binom{n}{2}}\frac{\poq{a}{n}}{\poq{q,c}{n}}\frac{nc^{n}}{a^{n+1}}\bigg|_{a=c}\\ &+\sum_{n=1}^\infty(-1)^nq^{\binom{n}{2}}\frac{\poq{a}{n}}{\poq{q,c}{n}}\bigg(\frac{c}{a}\bigg)^n\bigg\{\sum_{k=0}^{n-1}\frac{-q^k}{1-aq^k}\bigg\}\bigg|_{a=c}
\\
&=\bigg(\frac{c(c q/a; q)_{\infty}}{a^{2}(c ; q)_{\infty}}+(1-c/a)\times\mbox{other terms}\bigg)\bigg|_{a=c},
\end{align*}
which turns out to be
 \begin{align*}
-\sum_{n=1}^\infty(-1)^nq^{\binom{n}{2}}\frac{n}{\poq{q}{n}}+\sum_{n=1}^\infty(-1)^nq^{\binom{n}{2}}\frac{1}{\poq{q}{n}}\{\sum_{k=0}^{n-1}\frac{-cq^k}{1-cq^k}\}
=\frac{(q; q)_{\infty}}{(c ; q)_{\infty}}.
\end{align*}
A bit series rearrangement yields
 \begin{align*}-
\sum_{n=1}^\infty(-1)^nq^{\binom{n}{2}}\frac{1}{\poq{q}{n}}\{\sum_{k=0}^{n-1}\frac{1}{1-cq^k}\}
=\frac{(q; q)_{\infty}}{(c ; q)_{\infty}}.
\end{align*}
As desired.
\qed
\subsection{$q$-Identities based on the Gauss $_2\phi_1$ formula}
The famous $q$-Gauss formula states
\begin{lemma}[ \url{http://dlmf.nist.gov/17.6.E1}] For complex numbers $a,b,c: |c/ab|<1$, it holds
\begin{align}
{ }_2 \phi_1\left[\begin{array}{c}
a,b\\
c
\end{array} ;q, \frac{c}{ab}\right]=\frac{(c / a,c/b; q)_{\infty}}{(c,c/ab ; q)_{\infty}}.\label{wangjinyijian}
\end{align}
\end{lemma}By making use of this basic formula, we are able to set up
\begin{cor}
For complex number  $|b|>1$, it holds
\begin{align}\sum_{n=1}^\infty\frac{\poq{b}{n}}{\poq{q}{n}}\bigg(\frac{1}{b}\bigg)^n\bigg\{\sum_{k=0}^{n-1}\frac{1}{1-cq^k}\bigg\}
=-\frac{(q,c/b; q)_{\infty}}{(c,1/b; q)_{\infty}}.
\end{align}
\end{cor}
\pf~It suffices to  differentiate both sides of \eqref{wangjinyijian} with respect to  $a$ at $a=c.$ Then we have
 \begin{align*}
-\sum_{n=1}^\infty\frac{\poq{a,b}{n}}{\poq{q,c}{n}}n\frac{c^{n}}{a^{n+1}b^n}\bigg|_{a=c}&+\sum_{n=1}^\infty\frac{\poq{b}{n}}{\poq{q}{n}}\bigg(\frac{1}{b}\bigg)^n\{\sum_{k=0}^{n-1}\frac{-q^k}{1-cq^k}\}
=\frac{(q,c/b; q)_{\infty}}{c(c,1/b; q)_{\infty}}.
\end{align*}
That is
 \begin{align*}
\sum_{n=1}^\infty\frac{\poq{b}{n}}{\poq{q}{n}}\bigg(\frac{1}{b}\bigg)^n\bigg\{\sum_{k=0}^{n-1}\frac{1}{1-cq^k}\bigg\}
=-\frac{(q,c/b; q)_{\infty}}{(c,1/b; q)_{\infty}}.
\end{align*}
\qed
\begin{cor}
For complex number  $|c/b|<1$, it holds
\begin{align}\sum_{n=1}^{\infty}\frac{1}{1-q^n}\frac{\poq{b}{n}}{\poq{c}{n}}\left(\frac{c}{b}\right)^n=
\frac{\poq{c/b}{\infty}}{\poq{c}{\infty}}\sum_{n=0}^{\infty}n\frac{\poq{b}{n}}{\poq{q}{n}} \left(\frac{c}{b}\right)^{n}.
\end{align}
\end{cor}
\pf~ Actually,
we can differentiate both sides of \eqref{wangjinyijian} with respect to  $a$ at $a=1.$ Then we have
\begin{align*}
-\sum_{n=1}^{\infty}\frac{\poq{b}{n}}{\poq{c}{n}}\frac{\poq{aq}{n-1}}{\poq{q}{n}}\left(\frac{c}{ab}\right)^n\bigg|_{a=1}=
-\frac{\poq{c/b}{\infty}}{\poq{c}{\infty}}\sum_{n=0}^{\infty}\frac{\poq{b}{n}}{\poq{q}{n}}n \left(\frac{c}{ab}\right)^{n} \frac{1}{a}\bigg|_{a=1}.
\end{align*}
Therefore
\begin{align*}
\sum_{n=1}^{\infty}\frac{1}{1-q^n}\frac{\poq{b}{n}}{\poq{c}{n}}\left(\frac{c}{b}\right)^n=
\frac{\poq{c/b}{\infty}}{\poq{c}{\infty}}\sum_{n=0}^{\infty}\frac{\poq{b}{n}}{\poq{q}{n}}n \left(\frac{c}{b}\right)^{n}.
\end{align*}
This is what we are looking for.
\qed
\subsection{$q$-Identities based on Ramanujan's ${}_1\psi_1$ formula}
 Ramanujan's famous ${}_1\psi_1$  formula  can be recorded as follows.
\begin{lemma}[cf. \url{https://dlmf.nist.gov/17.8.E2}]For any complex numbers $a,b,x: abx\neq 0$, it holds
\begin{align}
{ }_1 \psi_1\left[\begin{array}{l}
a \\
b
\end{array} ; q, x\right]=\frac{(q, b / a, a x, q /a x; q)_{\infty}}{(b, q / a, x, b /ax; q)_{\infty}}
 \label{ramanujan}.
\end{align}
\end{lemma}

Performing as before, the first-order derivative of \eqref{ramanujan} leads us to
\begin{cor}For $|a|,|b|<1$, we have
\begin{align}
\frac{(q;q)_{\infty}^3(a b;q)_{\infty}}{(a;q)_{\infty}^2(b;q)_{\infty}^2}=\sum_{n=0}^{\infty}\left\{\frac{n(q / a;q)_n a^n}{(b;q)_{n+1}}+\frac{(n+1)(q / b;q)_n b^n}{(a;q)_{n+1}}\right\}.\label{danbian-gen-gen}
\end{align}
\end{cor}
\pf~  As carried out as above,  we apply $\mathbf{D}_{x=q/a}$  to both sides of \eqref{ramanujan}.
Then we obtain
\begin{align*}
\sum_{n=-\infty}^{\infty}\frac{(a;q)_n}{(b;q)_n}n \left(\frac{q}{a}\right)^{n-1}=\frac{a}{q}\frac{( b / a;q)_{\infty}}{(b, b/q;q)_{\infty}}\frac{(q; q)_{\infty}^3}{(q/a ; q)^2_{\infty}}.
\end{align*}
Take $(a,b) \to (q/a, bq)$. Then we have
\begin{align*}
\sum_{n=-\infty}^{\infty}\frac{(q/a;q)_n}{(bq;q)_n}n a^{n}= \frac{(ab;q)_{\infty}}{(bq,b;q)_{\infty}}\frac{(q; q)_{\infty}^3}{(a ; q)^2_{\infty}},
\end{align*}
namely
\begin{align*}
 \sum_{n=0}^{\infty}\frac{(q/a;q)_n}{(bq;q)_n}n a^{n}+\sum_{n=1}^{\infty}\frac{(1/b;q)_n}{(a;q)_n}b^n (-n)= \frac{(ab;q)_{\infty}}{(bq,b;q)_{\infty}}\frac{(q; q)_{\infty}^3}{(a ; q)^2_{\infty}}.
\end{align*}
By dividing both sides with $1-b$, we obtain
\begin{align*}
 \sum_{n=0}^{\infty}\frac{(q/a;q)_n}{(b;q)_{n+1}}n a^{n}+\sum_{n=1}^{\infty}\frac{(q/b;q)_{n-1}}{(a;q)_n}n b^{n-1} = \frac{(ab;q)_{\infty}}{(b;q)^2_{\infty}}\frac{(q; q)_{\infty}^3}{(a ; q)^2_{\infty}}.
\end{align*}
It turns out to be \eqref{danbian-gen-gen}.
\qed

Identity \eqref{danbian-gen-gen} contains the following special cases which are of interest.
\begin{exam}Evidently, for $a=b$ in \eqref{danbian-gen-gen}, we have
\begin{align}\sum_{n=0}^{\infty}(2n+1)\frac{(q / a;q)_n a^n}{(a;q)_{n+1}}=
\frac{(q;q)_{\infty}^3(a^2;q)_{\infty}}{(a;q)_{\infty}^4}.
\end{align}
Further, letting $a\to 0$, we get the following well-known $q$-identity
\begin{align*}
\sum_{n=0}^{\infty}(2n+1)(-1)^n  q^{n+1 \choose 2} =(q; q)_{\infty}^3.
\end{align*}
Alternately,  take  $a,b $ such that $ab=q$ in \eqref{danbian-gen-gen}. Then  it is clear  that for  $0<|q|<|b|<1,$  \begin{align}\sum_{n=-\infty}^{\infty}\frac{nb^n}{1-q^{n+1}/b}=
\frac{(q;q)_{\infty}^4}{(b;q)_{\infty}^2(q/b;q)_{\infty}^2}.\label{danbian-gen}
\end{align}
In \eqref{danbian-gen}, replace $q$ with $q^2$ and let $b=-q$. Then we have
\begin{align}\sum_{n=-\infty}^{\infty}\frac{n(-q)^n}{1+q^{2n+1}}=
(q;q^2)_{\infty}^4(q^4;q^4)_{\infty}^4.\label{danbian-gen-000}
\end{align}
\end{exam}
\subsection{$q$-Identities based on the quintuple product identity}As one of basic results in $q$-series theory,  the  quintuple product identity
can  date back to a century ago. We refer the reader to S. Cooper's wonderful survey \cite{cooper} for its history and various proofs.
\begin{lemma}[The quintuple product identity: \url{https://dlmf.nist.gov/17.8.E3}]\begin{align}\sum_{n=-\infty}^{\infty} q^{3{n\choose 2}}(1-x q^n)\left(q x^3\right)^n=
(q, x, q / x ; q)_{\infty}\left(q x^2, q /x^2 ; q^2\right)_{\infty}.\label{quintuplee}
\end{align}
\end{lemma}

\begin{cor}We have\begin{align}6\sum_{n=-\infty}^{\infty} nq^{3{n\choose 2}+2n}=(q ; q)_{\infty}^3 \left(q ; q^2\right)_{\infty}^2-(-q,-q^2,q^3; q^3)_{\infty}.\label{216}
\end{align}
\end{cor}
\pf~ It suffices to apply $\mathbf{D}_x$ to \eqref{quintuplee}. The result is
\begin{align*}
\mathbf{D}_x(\mbox {LHS of \eqref{quintuplee}})
&=\sum_{n=-\infty}^{\infty} q^{3{n\choose 2}}\left(- q^{2n}x^{3n}+3n(1-x q^n)q^nx^{3n-1}\right);\\
\mathbf{D}_x(\mbox {RHS of \eqref{quintuplee}}) &=(q ; q)_{\infty} \mathbf{D}_x((1-x)  F(x))\\
 &=(q ; q)_{\infty} \left\{-F(x)+(1-x)  \mathbf{D}_x(F(x))\right\},
 \end{align*}
 where $F(x)=(xq, q / x ; q)_{\infty}\left(q x^2, q /x^2 ; q^2\right)_{\infty}.$
Next, consider each derivative at $x=1$.
Consequently,
\begin{align*}
\mathbf{D}_{x=1}(\mbox {LHS of \eqref{quintuplee}})=\sum_{n=-\infty}^{\infty}\left(-q^{2n}+3 n(1- q^n)q^{n}\right) q^{3{n\choose 2}}\qquad\mbox{and}\\
\mathbf{D}_{x=1}(\mbox {RHS of \eqref{quintuplee}})=-(q ; q)_{\infty} F(1)= -(q ; q)_{\infty}^3 \left(q ; q^2\right)_{\infty}^2.
\end{align*}
Therefore
\begin{align}\sum_{n=-\infty}^{\infty}\left(1-3 n(1- q^n)q^{-n}\right) q^{3{n\choose 2}+2n}= (q ; q)_{\infty}^3 \left(q ; q^2\right)_{\infty}^2.\label{tttt}
\end{align}Furthermore,  it is easy to compute  that
\begin{align*} \mbox{LHS of \eqref{tttt}}
=\sum_{n=-\infty}^{\infty} \left(1+3n\right) q^{3{n\choose 2}+2n}-\sum_{n=-\infty}^{\infty}3n q^{3{n\choose 2}+n}.
\end{align*}Since, by changing the index of summation $n$ to $-n$,
\[-\sum_{n=-\infty}^{\infty}3n q^{3{n\choose 2}+n}=\sum_{n=-\infty}^{\infty}3n q^{3{n\choose 2}+2n},\]
we have
\begin{align*} \sum_{n=-\infty}^{\infty}\left(1-3 n(1- q^n)q^{-n}\right) q^{3{n\choose 2}+2n}
=\sum_{n=-\infty}^{\infty}  q^{3{n\choose 2}+2n}+6\sum_{n=-\infty}^{\infty} n q^{3{n\choose 2}+2n}.
\end{align*}As a last step, by Jacobi's triple product identity (cf. \url{https://dlmf.nist.gov/17.8.E1})
\begin{align}
\sum_{n=-\infty}^\infty(-1)^nq^{\binom{n}{2}}x^n=\poq{x,q/x,q}{\infty},\label{shuangbian}
\end{align}
 it is easy to find  $$\sum_{n=-\infty}^{\infty}  q^{3{n\choose 2}+2n}=(-q,-q^2,q^3;q^3)_\infty$$
Substituting this in \eqref{tttt}, we have \eqref{216} as desired.
\qed

\section{Main results-II: $q$-Identities from some advanced ${}_r\phi_s$ formulas}
In this part, we will proceed to consider some  advanced  ${}_r\phi_s$ formulas from \url{https://dlmf.nist.gov/17}. The word ``advanced'' means that these formulas often contain many parameters and their proofs  seem more complicated than those in the previous section.
\subsection{$q$-Identities based on the $q$-Clausen product formula}
The $q$-Clausen product formula due to G. Gasper can be recorded as
\begin{lemma}[cf. \url{https://dlmf.nist.gov/17.9.E18}]\label{aaaa}
\begin{align}\left({ }_4 \phi_3\left[\begin{array}{c}
a, b, a b x, a b / x \\
a b q^{\frac{1}{2}},-a b q^{\frac{1}{2}},-a b
\end{array}; q, q\right]\right)^2={ }_5 \phi_4\left[\begin{array}{c}
a^2, b^2, a b, a b x, a b / x \\
a b q^{\frac{1}{2}},-a b q^{\frac{1}{2}},-a b, a^2 b^2
\end{array} ; q, q\right],\label{kkkk}
\end{align}
where   $a, b$ with at least one is of form  $q^{-n}$, $n$ being nonnegative integer.
\end{lemma}
Appealing  to this product formula, we can show
\begin{cor}
For complex numbers  $a, b$ with at least one is of form  $q^{-n},n\geq 1$, it holds
\begin{align}
2\sum_{k=1}^{n}\frac{\poq{a,b,a^2b^2}{k}}{\poq{-ab}{k}\poqq{a^2b^2q}{k}}\frac{q^k}{1-q^k}=
\sum_{k=1}^{2n}\frac{\poq{a^2,b^2,ab}{k}}{\poq{-ab}{k}\poqq{a^2b^2q}{k}}\frac{q^k}{1-q^k}.\label{212}
\end{align}
\end{cor}
\pf~Observe that \eqref{kkkk} can be restated as
\begin{align}\bigg(1+
\sum_{n=1}^{\infty}\frac{\poq{a,b,abx,ab/x}{n}}{\poq{q,-ab}{n}\poqq{a^2b^2q}{n}}q^n\bigg)^2=
\sum_{n=0}^{\infty}\frac{\poq{a^2,b^2,ab,abx,ab/x}{n}}{\poq{q,-ab,a^2b^2}{n}\poqq{a^2b^2q}{n}}q^n.\label{xxxx-new}
\end{align}
Next, differentiate both sides of \eqref{xxxx-new} with respect to  $x$ at $x=ab.$ Then we have
 \begin{align*}
2\bigg(1+&
\sum_{n=1}^{\infty}\frac{\poq{a,b,abx,ab/x}{n}}{\poq{q,-ab}{n}\poqq{a^2b^2q}{n}}q^n\bigg)\frac{ab}{x^2}\sum_{n=1}^{\infty}\frac{\poq{a,b,abx}{n}\poq{abq/x}{n-1}}{\poq{q,-ab}{n}\poqq{a^2b^2q}{n}}q^n\bigg|_{x=ab}\\
&=\frac{ab}{x^2}
\sum_{n=1}^{\infty}\frac{\poq{a^2,b^2,ab,abx}{n}\poq{abq/x}{n-1}}{\poq{q,-ab,a^2b^2}{n}\poqq{a^2b^2q}{n}}q^n\bigg|_{x=ab},
\end{align*}
which turns out to be
  \begin{align*}
2\sum_{n=1}^{\infty}\frac{\poq{a,b,a^2b^2}{n}\poq{q}{n-1}}{\poq{q,-ab}{n}\poqq{a^2b^2q}{n}}q^n=
\sum_{n=1}^{\infty}\frac{\poq{a^2,b^2,ab,a^2b^2}{n}\poq{q}{n-1}}{\poq{q,-ab,a^2b^2}{n}\poqq{a^2b^2q}{n}}q^n.
\end{align*}
That is
 \begin{align*}
2\sum_{n=1}^{\infty}\frac{\poq{a,b,a^2b^2}{n}}{\poq{-ab}{n}\poqq{a^2b^2q}{n}}\frac{q^n}{1-q^n}=
\sum_{n=1}^{\infty}\frac{\poq{a^2,b^2,ab}{n}}{\poq{-ab}{n}\poqq{a^2b^2q}{n}}\frac{q^n}{1-q^n}.
\end{align*}
This  reduces to \eqref{212} immediately after the assumption of Lemma \ref{aaaa} is taken into account.
\qed

\subsection{$q$-Identities based on Rogers' $_6\phi_5$ formula}
\begin{lemma}[cf.  \url{https://dlmf.nist.gov/17.7.E15}]For $|aq/(bcd)|<1$, it holds
\begin{align}
{ }_6 \phi_5\left[\begin{array}{c}
a, q a^{\frac{1}{2}},-q a^{\frac{1}{2}}, b, c, d \\
a^{\frac{1}{2}},-a^{\frac{1}{2}}, a q / b, a q / c, a q / d
\end{array} ; q, \frac{a q}{b c d}\right]=\frac{(a q, a q /b c, a q /b d, a q /c d; q)_{\infty}}{(a q / b, a q / c, a q / d, a q /(b c d) ; q)_{\infty}}.\label{rogers65}
\end{align}
\end{lemma}
This frequently used formula implies
\begin{cor}[cf. \url{https://dlmf.nist.gov/17.8.E5}]
For complex number  $|q/(bcd)|<1$, it holds
\begin{align}{ }_3 \psi_3\left[\begin{array}{c}
b,c,d \\
q/b,q/c,q/d
\end{array}; q, \frac{q}{bcd}\right]=\frac{(q, q/b c,q/b d, q /c d; q)_{\infty}}{(q/ b,q/ c, q/ d, q/(b c d) ; q)_{\infty}};\label{PPPPP}
\end{align}
\begin{align}
C_0&+\sum_{n=2}^{\infty}\frac{(1-q^{2n-1})}{(1-q^{n-1})(1-q^n)}\frac{\poq{b,c,d}{n}}{\poq{1/b,1/c,1/d}{n}}\bigg(\frac{1}{b c d}\bigg)^n\nonumber\\
&=\frac{(q, 1 /b c, 1 /b d, 1 /c d; q)_{\infty}}
{(1 / b, 1/ c, 1 / d, 1 /(b c d) ; q)_{\infty}},\label{zzzz}
\end{align}
where the constant
$$C_0:=\frac{2 b c d q-b c d-b c q-b d q-c d q+b+c+d+q-2}{(1-q)(1-b) (1-c) (1-d)}.$$
In particular,  we have
\begin{align}&\sum_{n=2}^{\infty}
\frac{1-q^{2n-1}}{(1-q^{n-1})(1-q^n)}\frac{\poq{b}{n}^3}{\poq{1/b}{n}^3}\bigg(\frac{1}{b}\bigg)^{3n}\nonumber\\
&\qquad=(q;q)_\infty\frac{(-1 /b; q)_{\infty}^3(q/b^2;q^2)_\infty^3}{(1 /b^3; q)_{\infty}}+
\frac{-2 b q+b-q+2}{(1-b)(1-q)};\label{zzzz-old}
\\
&\sum_{n=2}^{\infty}\frac{1-q^{2n-1}}{(1-q^{n-1})(1-q^n)}(-1)^nq^{3n(n-1)/2}=(q; q)_{\infty}-\frac{1-2q}{1-q}.\label{zzzz-special}
\end{align}
\end{cor}
\pf~ Obviously, \eqref{rogers65} is equivalent to
\begin{align}1-a&+
\sum_{n=1}^{\infty}(1-aq^{2n})\frac{\poq{a,b,c,d}{n}}{\poq{q,aq/b,aq/c,aq/d}{n}}\bigg(\frac{a q}{b c d}\bigg)^n\nonumber\\
&=\frac{(a, a q /b c, a q /b d, a q /c d; q)_{\infty}}
{(a q / b, a q / c, a q / d, a q /(b c d) ; q)_{\infty}}.\label{xxxx-1}
\end{align}
Next, differentiate both sides of \eqref{xxxx-1} with respect to  $a$ at different points, say $a=1$ and $a=1/q.$
\begin{description}
\item[(i) at $a=1$]  in this case,  we have
 \begin{align*}
&-1-\sum_{n=1}^{\infty}(1-q^{2n})\frac{\poq{q}{n-1}}{\poq{q}{n}}\frac{\poq{b,c,d}{n}}{\poq{q/b,q/c,q/d}{n}}\bigg(\frac{q}{bcd}\bigg)^n
\\
&=-\frac{(q, q/b c,q/b d, q /c d; q)_{\infty}}{(q/ b,q/ c, q/ d, q/(b c d) ; q)_{\infty}},
\end{align*}
which turns out to be
  \begin{align*}
1+\sum_{n=1}^{\infty}(1+q^{n})\frac{\poq{b,c,d}{n}}{\poq{q/b,q/c,q/d}{n}}\bigg(\frac{q}{bcd}\bigg)^n
=\frac{(q, q/b c,q/b d, q /c d; q)_{\infty}}{(q/ b,q/ c, q/ d, q/(b c d) ; q)_{\infty}},
\end{align*}
which, written in terms of bilateral series, is \eqref{PPPPP}.
\item[(ii) at $a=1/q$] so we have
 \begin{align}
S(a)&+(1-a)(1-aq)
\sum_{n=2}^{\infty}(1-aq^{2n})\frac{\poq{aq^{2}}{n-2}\poq{b,c,d}{n}}{\poq{q,aq/b,aq/c,aq/d}{n}}\bigg(\frac{a q}{b c d}\bigg)^n\nonumber\\
&=\frac{(a, a q /b c, a q /b d, a q /c d ; q)_{\infty}}
{(a q / b, a q / c, a q / d, a q /(b c d) ; q)_{\infty}},\label{xxxx-xxx}
\end{align}
where
\[S(a):=1-a+aq(1-aq^2)\frac{(1-a)(1-b)(1-c)(1-d)}{(1-q)(b-aq)(c-aq)(d-aq)}.\]
Now differentiate both sides of \eqref{xxxx-xxx} with respect to  $a$ at  $a=1/q.$  The result is as follows:
\begin{align*}
&\frac{-2 b c d q+b c d+b c q+b d q-b+c d q-c-d-q+2}{(b-1) (c-1) (d-1)}\\
&-q(1-q^{-1})
\sum_{n=2}^{\infty}\frac{(1-q^{2n-1})}{(1-q^{n-1})(1-q^n)}\frac{\poq{b,c,d}{n}}{\poq{1/b,1/c,1/d}{n}}\bigg(\frac{1}{b c d}\bigg)^n\\
&=-q(1-q^{-1})\frac{(q, 1 /b c, 1 /b d, 1 /c d; q)_{\infty}}
{(1 / b, 1/ c, 1 / d, 1 /(b c d) ; q)_{\infty}},
\end{align*}
which turns out to be \eqref{zzzz}.
\end{description}
 We remark that \eqref{zzzz-old} and \eqref{zzzz-special} follow  from \eqref{zzzz}  with $b=c=d$ and $b,c,d\to \infty,$ respectively.\qed

\subsection{$q$-Identity based on Watson's  $q$-analog of Whipple's  theorem}
\begin{lemma}[cf. \url{http://dlmf.nist.gov/17.9.E15}]For  any integer $n\geq 0$, it holds
\begin{align}
& { }_8 \phi_7\left[
{{a, q a^{\frac{1}{2}},-q a^{\frac{1}{2}}, b, c, d, e, q^{-n}} \atop
{a^{\frac{1}{2}},-a^{\frac{1}{2}}, a q / b, a q / c, a q / d, a q / e, a q^{n+1}}} ; q, \frac{a^2 q^{2+n}}{b c d e}
\right]\label{whipplethm}\\
& =\frac{(a q, a q / d e ; q)_n}{(a q / d, a q / e ; q)_n}{ }_4 \phi_3\left[\begin{array}{c}
q^{-n}, d, e, a q / b c \\
a q / b, a q / c, d e q^{-n} / a
\end{array} ; q, q\right].\nonumber
\end{align}
\end{lemma}
From this fundamental result we can derive
\begin{cor}For any integer $n\geq 0$, it holds
\begin{align}
{}_5\psi_5\left[{{b, c, d, e, q^{-n}}\atop { q / b,  q / c,  q / d,  q / e,  q^{n+1}}};q,\frac{ q^{2+n}}{b c d e}\right]
=\frac{( q,  q / d e ; q)_n}{( q / d,  q / e ; q)_n}{ }_4 \phi_3\left[{{
q^{-n}, d, e,  q / b c} \atop
{ q / b,  q / c, d e q^{-n} }}
; q, q\right].\label{need}
\end{align}
In particular, when $bcdeq^{-n}=q$, it holds
\begin{align}
{ }_5 \psi_5\left[\begin{array}{c}b, c, d, e, q^{-n} \\
 q / b, q / c, q / d,  q / e,  q^{n+1}
\end{array}; q, q\right]=\frac{(q, q /b d, q /b e,  q /de; q)_n}{(q / b, q / d,  q / e,  q /(bde) ; q)_n}.\label{final-one}
\end{align}
\end{cor}
\pf~ To show \eqref{need}, we now reformulate \eqref{whipplethm} in the form
\begin{align}
& 1-a+\sum_{k=1}^{\infty} (1-a q^{2k})
\frac{\poq{a, b, c, d, e, q^{-n}}{k}} {\poq{q, a q / b, a q / c, a q / d, a q / e, a q^{n+1}}{k}} \left(\frac{a^2 q^{2+n}}{b c d e}\right)^k\nonumber\\
& =(1-a)\frac{(a q, a q / d e ; q)_n}{(a q / d, a q / e ; q)_n}{ }_4 \phi_3\left[\begin{array}{c}
q^{-n}, d, e, a q / b c \\
a q / b, a q / c, d e q^{-n} / a
\end{array} ; q, q\right].\label{wangjin}
\end{align}
Next, differentiate both sides of \eqref{wangjin} with respect to  $a$ at $a=1$. Then we have
\begin{align*}
\mathbf{D}_{a=1}(\mbox{LHS of  \eqref{wangjin}})&=-1-\sum_{k=1}^{\infty}(1+q^k) \frac{\poq{ b, c, d, e, q^{-n}}{k}} {\poq{  q / b,  q / c,  q / d,  q / e,  q^{n+1}}{k}}\left(\frac{ q^{2+n}}{b c d e}\right)^k\\
&=-\sum_{k=-n}^{n} \frac{\poq{ b, c, d, e, q^{-n}}{k}} {\poq{  q / b,  q / c,  q / d,  q / e,  q^{n+1}}{k}}\left(\frac{ q^{2+n}}{b c d e}\right)^k.
\end{align*}On the other hand,
\begin{align*}
\mathbf{D}_{a=1}(\mbox{RHS of \eqref{wangjin}})
=-\frac{( q,  q / d e ; q)_n}{( q / d,  q / e ; q)_n}{ }_4 \phi_3\left[{{
q^{-n}, d, e,  q / b c} \atop
{ q / b,  q / c, d e q^{-n} }}
; q, q\right],
\end{align*}
 we thus have
\begin{align*}
\sum_{k=-n}^{n} \frac{\poq{ b, c, d, e, q^{-n}}{k}} {\poq{  q / b,  q / c,  q / d,  q / e,  q^{n+1}}{k}}\left(\frac{ q^{2+n}}{b c d e}\right)^k
=\frac{( q,  q / d e ; q)_n}{( q / d,  q / e ; q)_n}{ }_4 \phi_3\left[{{
q^{-n}, d, e,  q / b c} \atop
{ q / b,  q / c, d e q^{-n} }}
; q, q\right].
\end{align*}
It amounts to \eqref{need}.

To  show  \eqref{final-one}, we see that for $bcdeq^{-n}=q$, the right-hand side of \eqref{need} becomes
\begin{align*}
\frac{( q,  q / d e ; q)_n}{( q / d,  q / e ; q)_n}{ }_3 \phi_2\left[{{
q^{-n}, d, e} \atop
{ q / b,  q / c}}
; q, q\right]=\frac{(q, q /b d, q /b e,  q /de; q)_n}{(q / b, q / d,  q / e,  q /(bde) ; q)_n}.
\end{align*}
The last equality comes from the $q$-Pfaff-Saalschutz ${}_3\phi_2$ formula (cf. \url{http://dlmf.nist.gov/17.7.E4}).
\qed

From the above transformation it follows a finite form of the famous Rogers-Ramanujan identity (see \eqref{ramanujanid} below).
\begin{cor1}For any integer $n\geq 0$, it holds
\begin{align}\sum_{k=-n}^{n}\left[2n \atop n+k\right]_q (-1)^k q^{(5k-1)k/2}=\left[2n\atop n\right]_q\sum_{k=0}^{n}\left[n \atop k\right]_q q^{k^2}.\label{eq29}
\end{align}
\end{cor1}
\pf~It only needs to put $b,c,d,e\to \infty$ in \eqref{need} by noting the fact
\begin{align*}
\lim_{x\to\infty}\frac{\poq{x}{k}}{x^k}=(-1)^k q^{k(k-1)/2}.
\end{align*}
Then it is easy to check that
\begin{align*}
\mbox{LHS of \eqref{need}}&=\sum_{k=-n}^{n}\frac{(q^{-n};q)_k}{(q^{n+1};q)_k} q^{2k^2+nk}\qquad\mbox{and}\\
\mbox{RHS of \eqref{need}}&=\poq{q}{n}\sum_{k=0}^{n}\frac{\poq{q^{-n}}{k}}{\poq{q}{k}}(-1)^k q^{k(k-1)/2+k+nk}.
\end{align*}After a bit simplification by  the basic relation
$$\frac{\poq{q^{-n}}{k}}{\poq{q}{k}} q^{nk}=(-1)^kq^{\binom{k}{2}}\left[n \atop k\right]_q ,$$
 we get
$$\sum_{k=-n}^{n}\left[n \atop k\right]_q\frac{(q;q)_k}{(q^{n+1};q)_k}  (-1)^k q^{(5k-1)k/2}=\sum_{k=0}^{n}\left[n \atop k\right]_q q^{k^2},$$
which coincides with  \eqref{eq29} .
\qed
\begin{exam}It is interesting to  see that
by putting $n\to \infty$ on both sides of \eqref{eq29}, we have $$\sum_{k=-\infty}^{\infty} (-1)^k q^{(5k-1)k/2}=\poq{q}{\infty}\sum_{k=0}^{\infty}\frac{ q^{k^2}}{\poq{q}{k}}.$$
On applying the Jacobi triple product identity \eqref{shuangbian}, we recover the Rogers-Ramanujan identity
$$(q^2,q^3, q^5;q^5)_{\infty}=\poq{q}{\infty}\sum_{k=0}^{\infty}\frac{ q^{k^2}}{\poq{q}{k}},$$
i.e.,
\begin{align}\sum_{k=0}^{\infty}\frac{q^{k^2}}{\poq{q}{k}}=\frac{(q^2,q^3, q^5;q^5)_{\infty}}{\poq{q}{\infty}}=\frac{1}{(q,q^4;q^5)_{\infty}}.\label{ramanujanid}
\end{align}In this sense,  \eqref{eq29} does offer a finite form of \eqref{ramanujanid}.
\end{exam}

\subsection{$q$-Identity based on Bailey's nonterminating extension of Jackson's  $_8\phi_7$ formula}
\begin{lemma}[cf. \url{https://dlmf.nist.gov/17.7.E17}]For  $a^2q=bcdef$, it holds
\begin{align}
{ }_8 \phi_7\left[\begin{array}{c}
a, qa^{\frac{1}{2}},-q a^{\frac{1}{2}}, b, c, d, e, f \\
a{ }^{\frac{1}{2}},-a^{\frac{1}{2}}, a q / b, a q / c, a q / d, a q / e, a q / f
\end{array};q,q\right] \label{wwww-1}\\
-\frac{b}{a} \frac{(a q, c, d, e, f, b q / a, b q / c, b q / d, b q / e, b q / f ; q)_{\infty}}{\left(a q / b, a q / c, a q / d, a q / e, a q / f, b c / a, b d / a, b e / a, b f / a, b^2 q / a ; q\right)_{\infty}}\nonumber \\
\times{ }_8 \phi_7\left[\begin{array}{c}
b^2 / a, q b a^{-\frac{1}{2}},-q b a^{-\frac{1}{2}}, b, b c / a, b d / a, b e / a, b f / a \\
b a^{-\frac{1}{2}},-b a^{-\frac{1}{2}}, b q / a, b q / c, b q / d, b q / e, b q / f
\end{array};q,q\right]\nonumber \\
=\frac{(a q, b / a, a q /(c d), a q /(c e), a q /(c f), a q /(d e), a q /(d f), a q /(e f) ; q)_{\infty}}{(a q / c, a q / d, a q / e, a q / f, b c / a, b d / a, b e / a, b f / a ;q)_{\infty}}.\nonumber
\end{align}
\end{lemma}\
By applying the OD-operator $\mathbf{D}_a$ to \eqref{wwww-1}, we easily deduce \begin{cor}
For complex numbers  $b,c,d,e,f: bcdef=q$, it holds
\begin{align}&_5\psi_5 \left[{{b,c,d,e,f}\atop {q/b,q/c,q/d,q/e,q/f}};q,q\right]\label{wwww}\\
&-\frac{b}{1+b}\frac{(q, bq,c, d, e, f, b q / c, b q / d, b q / e, b q / f ; q)_{\infty}}{\left(q / b,  q / c, q / d,  q / e, q / f, b c, b d, b e, b f , b^2q; q\right)_{\infty}}\nonumber \\
&\times \sum_{n=0}^{\infty}(1+bq^{n})\frac{\poq{b^2,b c , b d, b e, b f }{n}}{\poq{q,b q / c, b q / d, b q / e, b q / f}{n}}q^n\nonumber\\
&=\frac{(q, b,  q /c d, q /c e,  q /c f,  q /d e, q /d f, q /e f ; q)_{\infty}}{( q / c, q / d, q / e, q / f, b c, b d, b e , b f ;q)_{\infty}}.\nonumber
\end{align}
\end{cor}
\pf~ First of all, observe that  \eqref{wwww-1} can be written in full form
\begin{align}1-a+
&\sum_{n=1}^{\infty}(1-aq^{2n})\frac{\poq{a,b,c,d,e,f}{n}}{\poq{q,aq/b,aq/c,aq/d,aq/e,aq/f}{n}}q^n\nonumber\\
&-\frac{b}{a} \frac{(a, c, d, e, f, b q / a, b q / c, b q / d, b q / e, b q / f ; q)_{\infty}}{\left(a q / b, a q / c, a q / d, a q / e, a q / f, b c / a, b d / a, b e / a, b f / a, b^2 / a ; q\right)_{\infty}}\nonumber \\
&\times\sum_{n=0}^{\infty}(1-b^2q^{2n}/a)\frac{\poq{b^2/a,b, b c / a, b d / a, b e / a, b f / a }{n}}{\poq{q,b q / a, b q / c, b q / d, b q / e, b q / f}{n}}q^n\nonumber\\
&=\frac{(a, b / a, a q /c d, a q /c e, a q /c f, a q /d e, a q /d f, a q /e f ; q)_{\infty}}{(a q / c, a q / d, a q / e, a q / f, b c / a, b d / a, b e / a, b f / a ;q)_{\infty}}.
\label{yyyy}
\end{align}
Next, differentiate both sides of \eqref{yyyy} with respect to  $a$ at $a=1.$  By simplifying the resulted,  we obtain
   \begin{align*}
&1+\sum_{n\geq 1}(1+q^{n})\frac{\poq{b,c,d,e,f}{n}}{\poq{q/b,q/c,q/d,q/e,q/f}{n}}q^n
\\
&-b(1-b)\frac{(q, c, d, e, f, b q , b q / c, b q / d, b q / e, b q / f ; q)_{\infty}}{\left(q / b,  q / c, q / d,  q / e, q / f, b c, b d, b e, b f, b^2; q\right)_{\infty}}\nonumber \\
&\times\sum_{n=0}^{\infty}(1+bq^{n})\frac{\poq{b^2,b c , b d, b e, b f }{n}}{\poq{q,b q / c, b q / d, b q / e, b q / f}{n}}q^n\nonumber\\
&=\frac{(q, b,  q /c d, q /c e,  q /c f,  q /d e, q /d f, q /e f ; q)_{\infty}}{( q / c, q / d, q / e, q / f, b c, b d, b e , b f ;q)_{\infty}}.
\end{align*}
It is in agreement with \eqref{wwww}.
\qed

\begin{cor1}
For complex numbers  $b,d,e,f:~def=1, b\neq 0$, it holds
\begin{align}\sum_{n=0}^{\infty}\frac{(1+bq^{n})\poq{b d, b e, b f }{n}}{\poq{ b q / d, b q / e, b q / f}{n}}q^n=\frac{( b d, b e , b f ;q)_{\infty}-(b/d, b/e, b/f; q)_{\infty}}{b(1-d)(1- e)(1-f)(b q / d, b q / e, b q / f ; q)_{\infty} }.\label{BCA}
\end{align}
\end{cor1}
\pf~ It suffices to set   $ bc=q,~def=1$ in \eqref{wwww}, obtaining
\begin{align}&_3\psi_3 \left[{{d,e,f}\atop {q/d,q/e,q/f}};q,q\right]\\
&-\frac{b}{1+b}\frac{(q, bq,c, d, e, f, b q / c, b q / d, b q / e, b q / f ; q)_{\infty}}{\left(q / b,  q / c, q / d,  q / e, q / f,q, b d, b e, b f , b^2q; q\right)_{\infty}}F(b,d,e)\nonumber\\
&=\frac{( b,   q /c d, q /c e,  q /c f,  q /d e, q /d f, q /e f ; q)_{\infty}}{( q / c, q / d, q / e, q / f,  b d, b e , b f ;q)_{\infty}},\nonumber
\end{align}
where
\begin{align*}
F(b,d,e):=\sum_{n=0}^{\infty}(1+bq^{n})\frac{\poq{b d, b e, b f }{n}}{\poq{ b q / d, b q / e, b q / f}{n}}q^n.
\end{align*}
 By applying \eqref{PPPPP} to the left-hand side, we have
\begin{align*}\frac{( q/de,q/df, q /ef ; q)_{\infty}}{(q/d,q/e, q/f ; q)_{\infty}}
&-\frac{b}{1+b}\frac{( bq,c, d, e, f, b q / c, b q / d, b q / e, b q / f ; q)_{\infty}}{\left(q / b,  q / c, q / d,  q / e, q / f, b d, b e, b f , b^2q; q\right)_{\infty}}F(b,d,e)\nonumber\\
&=\frac{( b,   q /c d, q /c e,  q /c f,  q /d e, q /d f, q /e f  ; q)_{\infty}}{( q / c, q / d, q / e, q / f,  b d, b e , b f ;q)_{\infty}}.\nonumber
\end{align*}
Taking $bc=q,~def=1$ into account and simplifying further, we get
\begin{align}(dq,eq,fq; q)_{\infty}
&-b\frac{( d, e, f, b q / c, b q / d, b q / e, b q / f ; q)_{\infty}}{\left( b d, b e, b f , b^2; q\right)_{\infty}}F(b,d,e)\nonumber\\
&=\frac{( q /c d, q /c e,  q /c f,  q /d e, q /d f, q /e f  ; q)_{\infty}}{( b d, b e , b f ;q)_{\infty}}.\label{ABC}
\end{align}
Finally we solve  \eqref{ABC} for $F(b,d,e)$. The result is
\begin{align*}F(b,d,e)=\frac{( b d, b e , b f ;q)_{\infty}-(b/d, b/e, b/f; q)_{\infty}}{b(1-d)(1- e)(1-f)(b q / d, b q / e, b q / f ; q)_{\infty} }.\nonumber
\end{align*}
It completes the proof of \eqref{BCA}.
\qed

\subsection{Two further applications}
Now we turn to the following  double $q$-series transformation formula  due to Z. G. Liu \cite{30}. One of his proofs utilizes the $q$-exponential derivative operator. In \cite{Wang}, we have reproved it via the $(1-xy,y-x)$-expansion formula.
\begin{lemma}[{\rm cf. \cite[Thm. 1.5]{30} and \cite[Ex. 2.1]{Wang}}]Assume that $ |abcd/q^{2}|<1$. Then it holds
\begin{align} & \sum_{n=0}^{\infty} \frac{1-a q^{2 n}}{1-a} \frac{(a, q / b, q / c, q / d ; q)_n}{(q, a b, a c, a d ; q)_n}{ }_4 \phi_3\left[\begin{array}{c}q^{-n}, a q^n, \beta, \gamma \\ q / b, q / c, \beta \gamma a b c / q\end{array} ; q, q\right]\left(\frac{a b c d}{q^2}\right)^n\nonumber \\ & =\frac{\left(a q, \beta a b c / q, \gamma a b c / q, a b d / q, a c d / q, \beta \gamma a b c d / q^2 ; q\right)_{\infty}}{\left(a b, a c, a d, \beta \gamma a b c / q, \beta a b c d / q^2, \gamma a b c d / q^2 ; q\right)_{\infty}}.\label{KKK-1}\end{align}
\end{lemma}
Starting with this  transformation, we are now able to extend it to bilateral series.

\begin{cor}For  $|bcd/q^2|<1$, it holds
\begin{align} & \sum_{n=-\infty}^{\infty}G(n;\beta,\gamma) \frac{(q / b, q / c, q / d ; q)_n}{( b, c, d; q)_n}\left(\frac{ b c d}{q^2}\right)^n\nonumber \\ & =\frac{\left(q, \beta b c / q, \gamma b c / q, b d / q, c d / q, \beta \gamma b c d / q^2 ; q\right)_{\infty}}{\left( b, c, d, \beta \gamma b c / q, \beta b c d / q^2, \gamma b c d / q^2 ; q\right)_{\infty}},\label{KKK}\end{align}
where
\[G(n;\beta,\gamma):={ }_4 \phi_3\left[\begin{array}{c}q^{-n}, q^n, \beta, \gamma \\ q / b, q / c, \beta \gamma b c / q\end{array} ; q, q\right].\]
\end{cor}
\pf It suffices to multiply both sides of \eqref{KKK-1} with $1-a$, and then apply  $\mathbf{D}_{a=1}.$  Finally,  reformulate the left-hand infinite series of  the resulting identity in terms of bilateral  hypergeometric series. Note that $G(n;\beta,\gamma)=G(-n;\beta,\gamma)$.
\qed
\begin{remark} Note that  the case $\beta=1$ or $\gamma=1$ of  \eqref{KKK}  with the replacement  $(b,c,d)\to (q/b,q/c,q/d)$ reduces to \eqref{PPPPP}.
\end{remark}

At the end of this paper,  we should point out that the foregoing argument is also  applicable  for any proper  transformations   among   ${}_r\phi_s$  series. The following is a good example for this point.
\begin{lemma}[cf. \url{https://dlmf.nist.gov/17.9.E3_5}] For $|x|<1,$ it holds
\begin{align}
{ }_2 \phi_1\left[\begin{array}{c}
a, b \\
c
\end{array} ; q, x\right]&=\frac{(c / a, c / b ; q)_{\infty}}{(c, c /(a b) ; q)_{\infty}} {}_3\phi_2\left[\begin{array}{c}
a, b, a bx/ c \\
q a b / c, 0
\end{array} ; q, q\right]\label{QQQQ} \\
&+\frac{(a, b, a bx / c ; q)_{\infty}}{(c, a b / c, x; q)_{\infty}} {}_3\phi_2\left[\begin{array}{c}
c / a, c / b, x \\
q c /a b, 0
\end{array} ; q, q\right].\nonumber
\end{align}
\end{lemma}
By virtue of \eqref{QQQQ}, it is easy to establish
\begin{cor} For $|c/ab|<1$, it holds
\begin{align}
\sum_{n=0}^{\infty}n\frac{\poq{a,b}{n}}{\poq{q,c}{n}}\left(\frac{c}{ab}\right)^n
&=-
\frac{(c / a, c / b ; q)_{\infty}}{(c, c /(a b) ; q)_{\infty}}\sum_{n=1}^{\infty}\frac{\poq{a,b}{n}}{\poq{qab/c}{n}}\frac{q^{n}}{1-q^{n}}\label{QQQQ-new}\\
&-\frac{(a, b, q ; q)_{\infty}}{(c, a b / c, cq/ab ; q)_{\infty}}\sum_{n=0}^{\infty}\frac{\poq{c/a,c/b}{n}}{\poq{q}{n}}\frac{q^{n}}{1-cq^{n}/ab}.\nonumber
\end{align}
\end{cor}
\pf~It suffices to differentiate both sides of \eqref{QQQQ} with respect to  $x$ at $x=c/ab.$ Then we have
\begin{align*}
\sum_{n=1}^{\infty}\frac{\poq{a,b}{n}}{\poq{q,c}{n}}nx^{n-1}\bigg|_{x=c/ab}&=-\frac{ab}{c}
\frac{(c / a, c / b ; q)_{\infty}}{(c, c /(a b) ; q)_{\infty}}\sum_{n=1}^{\infty}\frac{\poq{a,b}{n}}{\poq{qab/c}{n}}\frac{\poq{abxq/c}{n-1}}{\poq{q}{n}}q^{n}\bigg|_{x=c/ab}  \\
&-\frac{ab}{c}\frac{(a, b, a bxq / c ; q)_{\infty}}{(c, a b / c, x ; q)_{\infty}}\sum_{n=0}^{\infty}\frac{\poq{c/a,c/b,x}{n}}{\poq{q,cq/ab}{n}}q^{n}\bigg|_{x=c/ab}.
\end{align*}
After some routine computation, it becomes
\begin{align*}
\sum_{n=0}^{\infty}n\frac{\poq{a,b}{n}}{\poq{q,c}{n}}\left(\frac{c}{ab}\right)^n&=-
\frac{(c / a, c / b ; q)_{\infty}}{(c, c /(a b) ; q)_{\infty}}\sum_{n=1}^{\infty}\frac{\poq{a,b}{n}}{\poq{qab/c}{n}}\frac{q^{n}}{1-q^{n}} \\
&-\frac{(a, b, q ; q)_{\infty}}{(c, a b / c, c/ab ; q)_{\infty}}\sum_{n=0}^{\infty}\frac{\poq{c/a,c/b,c/ab}{n}}{\poq{q,cq/ab}{n}}q^{n}.
\end{align*}
It gives the complete proof of \eqref{QQQQ-new}.
\qed
\begin{exam}For the case $b=c$, \eqref{QQQQ-new} reduces to
\begin{align*}
\sum_{n=0}^{\infty}n\frac{\poq{a}{n}}{\poq{q}{n}}\left(\frac{1}{a}\right)^n=-\frac{(q ; q)_{\infty}}{(1/a; q)_{\infty}}.
\end{align*}
It is also  the special case of \eqref{eq2.8} with $m=0$.
\end{exam}
%
\section*{Acknowledgements}
This  work was supported by the National Natural Science Foundation of
China [Grant Nos.  12001492 and 11971341].

\end{document}